\documentclass[11pt]{article}
\usepackage{latexsym}
\usepackage[tbtags]{amsmath}
\usepackage{epsfig,amstext,amssymb,amsthm,latexsym}
\usepackage{amssymb,color}

\UseRawInputEncoding

\definecolor{c20}{rgb}{0.,0.7,0.}
\definecolor{c30}{rgb}{0.,0.,1.}
\definecolor{c40}{rgb}{1,0.1,0.7}
\definecolor{c50}{rgb}{1,0,0}

\date{}
\newtheorem{theorem}{Theorem}[section]

\makeatletter % `@' is now normal letter for TeX
\@addtoreset{equation}{section}
\makeatother % `@' is restored as a ``non-letter'' character
\begin{document}
%%%%%%%%%%%%%%%%%%%%%%%%%%%%%%%%%%%%%%%%%%%%%%%%%%%%%%%%%
%%%%%%%%%%%%%%%%%%%%%%%%%%%%%%%%%%%%%%%%%%%%%%%%%%%%%%%%%
\title{Kendall Correlation Coefficient for non-Identically Distributed Variables}

\author{Alexei Stepanov \thanks{\noindent  Education and Research Cluster ``Institute of High Technology",\  Immanuel Kant  Baltic Federal University, A.Nevskogo 14, Kaliningrad, 236041 Russia,  email: alexeistep45@mail.ru}  }

\maketitle
\begin{abstract} 
In the present paper, we    discuss for the first time the theoretical Kendall correlation coefficient for non-identical bivariate data. In the non-identical case, we first introduce a theoretical Kendall correlation coefficient $\tau_n$ and  show that the expected value of the rank Kendall correlation coefficient $\tilde{\tau}_n$ is equal to $\tau_n$. We then prove that $\tilde{\tau}_n$ converges in probability to  $\tau=\lim_{n\rightarrow\infty} \tau_n$.  These facts enable us to state that  $\tau_n$ is a correctly defined  theoretical Kendall correlation coefficient for the non-identical case. We also support our theoretical results by simulation experiments.
\end{abstract}



\section{Introduction} Let  $(X,Y), (X_1,Y_1),\ldots,(X_n,Y_n)$  be independent and identically distributed (iid) random vectors with absolutely continuous bivariate distribution function $F(x,y)=P(X\leq x, Y\leq y)$,   density function  $f(x,y)$ and marginal distribution  functions    $H(x)$ and $G(y)$.  There are three basic  measures of dependence rate: the Pearson $\rho$, the Spearman $\rho_S$ and the Kendall $\tau$ correlation coefficients. 

The rate of dependence between the variables $X$ and $Y$ is most often measured by the Pearson correlation coefficient
$$
\rho=\frac{E(X-EX)(Y-EY)}{\sigma_X\sigma_Y}.
$$
The random variable
$$
\tilde{\rho}_n =\frac{\sum_{i=1}^n (X_i-\bar{X})(Y_i-\bar{Y})} {\sqrt{\sum_{i=1}^n (X_i-\bar{X})^2 \sum_{i=1}^n (Y_i-\bar{Y})^2}}
$$
represents the sample Pearson correlation coefficient, for which it is known that $\tilde{\rho}_n\stackrel{p}{\rightarrow}\rho$; see, for example,  Fisher (1921). 

Let  $X_{1,n}\leq\ldots\leq X_{n,n}$ be the order statistics obtained from the sample   $X_1,\ldots,X_n$. For these order statistics, let us define their concomitants   $Y_{1,n},\ldots,Y_{n,n}$. Let $X_i=X_{j,n}$, then  $Y_{j,n}=Y_i$ be the concomitant of the order statistic  $X_{j,n}$.  The concept of concomitants is   discussed in the papers of   David and Galambos (1974), Bhattacharya (1974), Egorov and Nevzorov (1984),  David and Nagaraja (2003), Balakrishnan and Lai (2009), Bairamov and Stepanov (2010),   Balakrishnan and Stepanov (2015). See also the references therein. 

Let  $I_{ji}=I(Y_{j,n}\leq Y_{i,n}),\ i\not= j\in \{1,2,\ldots,n\}$, where $I$ is the indicator function. The random variable  
\begin{equation}\label{1.1}
\tilde{\tau} _n=\frac{4\sum_{i=2}^n\sum_{j=1}^{i-1}I_{ji}}{n(n-1)}-1
\end{equation}
is known as the rank Kendall correlation coefficient. The theoretical analogue of $\tilde{\tau}_n$ is presented by
$$
\tau=4E[F(X,Y)]-1=4\int_{\mathbb{R}^2}F(x,y)f(x,y)dxdy-1.
$$
It is known, see, for example, Balakrishnan and Lai (2009) or Stepanov (2025a), that $\tau=E\tilde{\tau}_n$. It is also shown in Stepanov(2025a) that $\tilde{\tau}_n\stackrel{p}{\rightarrow} \tau$ as $n\rightarrow \infty$. The random variable  
$$
\tilde{\rho}_{n,S} =1-\frac{6\sum_{i=1}^n\left(\sum_{j=1}^{n}I_{ji}-i\right)^2}{n^3-n}
$$ 
is known as the rank Spearman correlation coefficient. The theoretical analogue of $\tilde{\rho}_{n,S}$ is presented by
$$
\rho_S =12E[H(X)G(Y)]-3=12\int_{\mathbb{R}^2}H(x)G(y)f(x,y)dxdy-3.
$$ 
New theoretical and rank correlation coefficients were introduced and analyzed in Stepanov (2025b). The  theoretical correlation coefficient  $r$ was defined by 
\begin{eqnarray*}
r&=&6E[F(X,Y)-H(X)G(Y)]\\
&=& 6\int_{\mathbb{R}^2}[F(x,y)-H(x)G(y)]f(x,y)dxdy
\end{eqnarray*}
and the rank correlation coefficient $\tilde{r}_n$  was defined by $\tilde{r}_n= \frac{3\tilde{\tau}_n-\tilde{\rho}_{n,S}}{2}$. The rationale for introducing $r$ is as follows. Let us consider the dependence function $F(x,y)-H(x)G(y)$, which gives us the dependence value in every point $(x,y)$ of the distribution support. We adjust the definition of the dependence function  and redefine it as $6[F(x,y)-H(x)G(y)]$. We can view any correlation coefficient     as  some average value of  the dependence rate, or of the dependence function. Observe that the dependence function does not take into account the likelihood of obtaining points $(x,y)$. In the above definition of $r$, we  take this likelihood into account. Thus, the correlation coefficient  $r$ represents  the ``ideal" average of the dependence rate. Since
$E\tilde{r}_n\rightarrow r,\  \tilde{r}_n\stackrel{p}{\rightarrow} r\ (n\rightarrow \infty)$, see Stepanov (2025b), the coefficient $\tilde{r}_n$ can be viewed as the rank representative of $r$.

The  correlation coefficients $\tilde{\rho}_n,\ \tilde{\rho}_{n,S},\ \tilde{\tau}_n$ and $\tilde{r}_n$ were compared in Stepanov (2025b). It  was revealed  that the behavior of $\tilde{\rho}_n$ can be very different from the behavior of the rank correlation coefficients $\tilde{\rho}_{n,S},\ \tilde{\tau}_n,\ \tilde{r}_n$, which, in turn, behave in a similar way with each other. This  follows from the definitions of $\tilde{\rho}_n$ and $\tilde{\rho}_{n,S},\ \tilde{\tau}_n,\ \tilde{r}_n$. The coefficient $\tilde{\rho}_n$ measures the association between the variables $X$ and $Y$, whereas  the coefficients $\tilde{\rho}_{n,S},\ \tilde{\tau}_n$ and $\tilde{r}_n$ measure the association between the corresponding ranks. The latter, in particular, means that the coefficients  $\tau,\ \rho_S$ and   $r$ are  invariant under any increasing transformation. 

A question was raised:  which  correlation coefficient measures the dependence rate best?  Based on the simulation analysis of the variances of correlation coefficients, it was concluded in Stepanov (2025b) that  $\tilde{\rho}_n$ is preferable when the association between $X$ and $Y$ is "close``\ to $\pm 1$ and  $\tilde{\rho}_{n,S},\ \tilde{\tau}_n,\ \tilde{r}_n$ are preferable, otherwise. Of course, the range "close``\ is negotiable.  It was also shown  that amongst the rank coefficients    $\tilde{\tau}_n,\ \tilde{\rho}_{n,S}$ and $\tilde{r}_n$ the coefficient  $\tilde{r}_n$  has the smallest variances when the association between $X$ and $Y$ is not "close``\ to $\pm 1$. Some   disadvantages of $\rho/\tilde{\rho}_n$ were  also highlight in Stepanov (2025b). For example, the disadvantage of $\rho$ and, correspondingly, of $\tilde{\rho}_n$ is that $\rho$ cannot be defined  if  the  second moments do not exist. Observe that $\rho_S,\ \tau$ and $r$ exist for any absolutely continuous non-degenerate distribution function. Another disadvantage of $\tilde{\rho}_n$ is that   $\tilde{\rho}_n$ is very sensitive to contamination of outliers, while $\tilde{\rho}_{n,S},\ \tilde{\tau}_n$ and $\tilde{r}_n$ are rather robust, see, for example, Shevlyakov and Vilchevski (2002),  Xu et al. (2013) and Stepanov (2025b).

The coefficients $\tilde{\rho}_n,\ \tilde{\rho}_{n,S},\ \tilde{\tau}_n$ and $\tilde{r}_n$ continue to be well-defined in the non-identical case. They can be easily computed, but their asymptotic behaviors are  unpredictable due to the difficulties with defining their theoretical analogues $\rho_n,\ \rho_{n,S},\ \tau_n$ and $r_n$. Observe that  the theoretical correlation coefficients depend on $n$ in the non-identical case.  The goal of this paper is to  determine and study (for the first time) the theoretical Kendall correlation coefficient in the non-identical case.

Let  $(X_1,Y_1),\ldots,(X_n,Y_n)$  be independent   random vectors with absolutely continuous bivariate distribution functions $F_1(x,y),\ldots, F_n(x,y)$,   density functions $f_1(x,y),\ldots, f_n(x,y)$ and marginal distribution  functions    $H_1(x),\ldots, H_n(x)$ and $G_1(y),\ldots,G_n(y)$.  The relationship  between the samples $X_1,\ldots, X_n$ and $Y_1,\ldots, Y_n$  can be measured by  the Kendall correlation coefficient defined by 
\begin{eqnarray}\label{1.2}
\tau_n&=&\frac{4\sum_{i=1}^n\sum_{j=1, j\not = i}^n E[F_i(X_j,Y_j)]}{n(n-1)}-1\\
\nonumber &=&\frac{4\sum_{i=1}^n\sum_{j=1, j\not = i}^n \int_{\mathbb{R}^2}F_i(x,y)f_j(x,y)dxdy}{n(n-1)}-1\\
\nonumber &=& \frac{4\sum_{i=1}^n\sum_{j=1, j\not = i}^n P(X_j\leq X_i,\ Y_j\leq  Y_i)}{n(n-1)}-1.
\end{eqnarray}
The rank correlation coefficient $\tilde{\tau}_n$ in the non-identical case can  again be defined by (\ref{1.1}). It has the same form as in the iid case. Observe that in the iid case,  i.e., when $F_i=F$, (\ref{1.2}) reduces to $\tau_n=4E[F(X,Y)]-1=\tau$. Since the sample $(X_i,Y_i)\ (i=1,\ldots,n)$ consists of non-identically distributed vectors, the coefficients $\tilde{\tau}_n$ and $\tau_n$ allows us also to work with   iid sample containing  outliers.

The rest of the paper is organized as follows. In Section~2, we show that $E\tilde{\tau}_n=\tau_n$ and that there exists a $\tau\in[0,1]$, such that
$$
\tau_n\rightarrow \tau,\quad \tilde{\tau}_n\stackrel{p}{\rightarrow} \tau\quad (n\rightarrow \infty).
$$
These results enable us to consider $\tilde{\tau}_n$ and $\tau_n$ as a pair of  conjugate Kendal correlation coefficients (rank and theoretical) in the non-identical case. In Section~3, we support our theoretical results by simulation experiments.

\section{Results}

\begin{theorem}\label{theorem2.1} The expected value of $\tilde{\tau}_n$ has the form
$$
E\tilde{\tau}_n=\tau_n.
$$
\end{theorem}

\begin{gproof}{ } 
For finding $E\tilde{\tau}_n$, let us study  $E\left[\sum_{i=2}^n\sum_{j=1}^{i-1}I_{ji}\right]=\sum_{i=2}^n\sum_{j=1}^{i-1}P(Y_{j,n}\leq Y_{i,n})$. 

For $n=2$, it has the form  $E\left[\sum_{i=2}^2\sum_{j=1}^{i-1}I_{ji}\right]=P(Y_{1,2}\leq Y_{2,2})$. We have two options: either $X_1=X_{1,2},\ X_2=X_{2,2}$, or    $X_1=X_{2,2},\ X_2=X_{1,2}$. Since $X_{1,2}\leq X_{2,2}$, we obtain that
\begin{eqnarray*}
P(Y_{1,2}\leq Y_{2,2}) &=& P(X_{1,2}\leq X_{2,2},\ Y_{1,2}\leq Y_{2,2})\\
&=& P(X_1\leq X_2, Y_1\leq Y_2)+P(X_2\leq X_1, Y_2\leq Y_1)=p_{21}+p_{12},\\
\end{eqnarray*}
where $p_{ij}=\int_{\mathbb{R}^2}F_j(x,y)f_i(x,y)dxdy$. The result is shown for $n=2$. Let us consider the case  $n=3$. We have
$$
E\left[\sum_{i=2}^3\sum_{j=1}^{i-1}I_{ji}\right]=P(Y_{1,3}\leq Y_{2,3})+P(Y_{1,3}\leq Y_{3,3})+P(Y_{2,3}\leq Y_{3,3}).
$$
Let
$$
B^{(3)}_{1,2}=\{X_{1,3}\leq X_{2,3},\ Y_{1,3}\leq Y_{2,3}\},\quad B^{(3)}_{1,3}=\{ X_{1,3}\leq X_{3,3},\ Y_{1,3}\leq Y_{3,3}\},
$$
$$
B^{(3)}_{2,3}=\{X_{2,3}\leq X_{3,3},\ Y_{2,3}\leq Y_{3,3}\}.
$$
We can write that
\begin{eqnarray*}
B^{(3)}_{1,2}&=&\cup_{1\leq m\leq 3,\ 1\leq k\leq 3,\ k\not= m}\{X_{1,3}\leq X_{2,3},\ Y_{1,3}\leq Y_{2,3},\ X_{1,3}=X_m,\ X_{2,3}=X_k\}\\
&=&\{X_1\leq X_2\leq X_3,\ Y_1\leq Y_2\} \cup \{X_2\leq X_1\leq X_3,\ Y_2\leq Y_1\}\\
&\cup&\{X_1\leq X_3\leq X_2,\ Y_1\leq Y_3\} \cup \{X_2\leq X_3\leq X_1,\ Y_2\leq Y_3\}\\
&\cup&\{X_3\leq X_1\leq X_2,\ Y_3\leq Y_1\} \cup \{X_3\leq X_2\leq X_1,\ Y_3\leq Y_2\},
\end{eqnarray*}
\begin{eqnarray*}
B^{(3)}_{1,3}&=&\cup_{1\leq m\leq 3,\ 1\leq k\leq 3,\ k\not= m}\{X_{1,3}\leq X_{3,3},\ Y_{1,3}\leq Y_{3,3},\ X_{1,3}=X_m,\ X_{3,3}=X_k\}\\
&=&\{X_1\leq X_3\leq X_2,\ Y_1\leq Y_2\} \cup \{X_2\leq X_3\leq X_1,\ Y_2\leq Y_1\}\\
&\cup&\{X_1\leq X_2\leq X_3,\ Y_1\leq Y_3\} \cup \{X_3\leq X_2\leq X_1,\ Y_3\leq Y_1\}\\
&\cup&\{X_2\leq X_1\leq X_3,\ Y_2\leq Y_3\} \cup \{X_3\leq X_1\leq X_2,\ Y_3\leq Y_2\}
\end{eqnarray*}
and
\begin{eqnarray*}
B^{(3)}_{2,3}&=&\cup_{1\leq m\leq 3,\ 1\leq k\leq 3,\ k\not= m}\{X_{2,3}\leq X_{2,3},\ Y_{3,3}\leq Y_{2,3},\ X_{2,3}=X_m,\ X_{3,3}=X_k\}\\
&=&\{X_3\leq X_1\leq X_2,\ Y_1\leq Y_2\} \cup \{X_3\leq X_2\leq X_1,\ Y_2\leq Y_1\}\\
&\cup&\{X_2\leq X_1\leq X_3,\ Y_1\leq Y_3\} \cup \{X_2\leq X_3\leq X_1,\ Y_3\leq Y_1\}\\
&\cup&\{X_1\leq X_2\leq X_3,\ Y_2\leq Y_3\} \cup \{X_1\leq X_3\leq X_2,\ Y_3\leq Y_2\}.
\end{eqnarray*}
Let
$$
A^{(3)}_{1,2}=\{X_2\leq X_1,\ Y_2\leq Y_1\},\quad A^{(3)}_{2,1}=\{X_1\leq X_2,\ Y_1\leq Y_2\},\quad A^{(3)}_{1,3}=\{X_3\leq X_1,\ Y_3\leq Y_1\},
$$
$$
A^{(3)}_{3,1}=\{X_1\leq X_3,\ Y_1\leq Y_3\},\quad A^{(3)}_{2,3}=\{X_3\leq X_2,\ Y_3\leq Y_2\},\quad A^{(3)}_{3,2}=\{X_2\leq X_3,\ Y_2\leq Y_3\}.
$$
Then
$$
B^{(3)}_{1,2}\cup B^{(3)}_{1,3}\cup B^{(3)}_{2,3}=A^{(3)}_{1,2}\cup A^{(3)}_{1,3}\cup A^{(3)}_{2,3}\cup A^{(3)}_{2,1}\cup A^{(3)}_{3,1}\cup A^{(3)}_{3,2}.
$$
It follows that
\begin{eqnarray*}
E\left[\sum_{i=2}^3\sum_{j=1}^{i-1}I_{ji}\right]&=&P\left(B^{(3)}_{1,2}\cup B^{(3)}_{1,3}\cup B^{(3)}_{2,3}\right)\\
&=&p_{12}+p_{21}+p_{13}+p_{31}+p_{23}+p_{32}.
\end{eqnarray*}
The result is shown for $n=3$. The above  argument can be applied further   for any integer $n\geq 4$. Let 
$$
A^{(n)}_{k,m}=\{X_m\leq X_k,\ Y_m\leq Y_k\}\quad \mbox{and}\quad B^{(n)}_{j,i}=\{X_{j,n}\leq X_{i,n},\ Y_{j,n}\leq Y_{i,n}\}.
$$
Since
$$
\cup_{1\leq j<i\leq n}\ B^{(n)}_{j,i}=\cup_{1\leq  m\leq n,\ 1\leq k\leq n,\ m\not= k}\ A^{(n)}_{k,m},
$$
we get
$$
E\left[\sum_{i=2}^n\sum_{j=1}^{i-1}I_{ji}\right]=\sum_{k=1}^n\ \sum_{m=1,\ m\not= k}^n p_{mk}.
$$
\end{gproof}

\begin{theorem}\label{theorem2.2} There exists a $\tau\in[0,1]$, such that
$$
\tau_n\rightarrow \tau\in [0,1]\quad (n\rightarrow \infty).
$$
\end{theorem}

\begin{gproof}{ } 
We will show that the sequence $\tau_n$ is fundamental. Observe that 
$$
\sum_{k=1}^{n+1}\ \sum_{m=1,\ m\not= k}^{n+1} p_{mk}=\sum_{k=1}^n\ \sum_{m=1,\ m\not= k}^n p_{mk}+\sum_{k=1}^n p_{(n+1)k}+\sum_{m=1}^n p_{m(n+1)}.
$$
We have
\begin{eqnarray*}
\mid \tau_{n+1}-\tau_n\mid &=&\frac{4\mid (n-1)\sum_{k=1}^{n+1}\ \sum_{m=1,\ m\not= k}^{n+1} p_{mk}-(n+1)\sum_{k=1}^{n}\ \sum_{m=1,\ m\not= k}^{n} p_{mk}\mid }{(n-1)n(n+1)}\\
&=&\frac{4\mid (n-1)\sum_{k=1}^n p_{(n+1)k}+(n-1)\sum_{m=1}^n p_{m(n+1)}-2\sum_{k=1}^{n}\ \sum_{m=1,\ m\not= k}^{n} p_{mk}\mid }{(n-1)n(n+1)}\rightarrow 0.
\end{eqnarray*}
In the same way, one can show that for any fixed $k\geq 1$ the convergence $\mid \tau_{n+k}-\tau_n\mid \rightarrow 0\ (n\rightarrow \infty)$ holds. The result readily follows.
\end{gproof}

\begin{theorem}\label{theorem2.3}  Let $\tau$ be a limit mentioned in Theorem~\ref{theorem2.2}. Then 
$$
\tilde{\tau}_n\stackrel{p}{\rightarrow} \tau\quad (n\rightarrow \infty).
$$
\end{theorem}

\begin{gproof}{ } It follows from Chebyshev's inequality that for any $\varepsilon >0$
\begin{eqnarray*} 
P(\mid \tilde{\tau}_n-\tau\mid> \varepsilon)&\leq &\frac{E(\tilde{\tau}_n-\tau)^2}{\varepsilon^2}\leq \frac{E(\tilde{\tau}_n-\tau_n+\tau_n-\tau)^2}{\varepsilon^2}\\
&=& \frac{Var(\tilde{\tau}_n)+2(\tau_n-\tau)E(\tilde{\tau}_n-\tau_n)+(\tau_n-\tau)^2}{\varepsilon^2}=\frac{Var(\tilde{\tau}_n)+o_n(1)}{\varepsilon^2}.
\end{eqnarray*}
We have
\begin{eqnarray}\label{2.1}
Var(\tilde{\tau}_n)\leq \frac{16}{(n-1)^4}\left(E\left[\sum_{i=2}^n\sum_{j=1}^{i-1}I_{ji}\right]^2-\left(\sum_{k=1}^n\ \sum_{m=1,\ m\not= k}^n p_{mk}\right)^2\right).
\end{eqnarray}
Observe that
$$
E\left[\sum_{i=2}^n\sum_{j=1}^{i-1}I_{ji}\right]^2=E\left[\sum_{i=2}^n\sum_{j=1}^{i-1}\sum_{k=2}^n\sum_{m=1}^{k-1}I_{ji}I_{mk}\right]=E_1+E_2+E_3+E_4+E_5,
$$
where 
$$
E_1=E\left[\sum_{i=2}^n\sum_{j=1}^{i-1}I^2_{ji}\right]\leq \frac{n(n-1)}{2}\quad (\mbox{here}\ m=j,\ k=i),
$$
\begin{eqnarray*}
E_2=2E\left[\sum_{k=3}^n\sum_{i=2}^{k-1}\sum_{j=1}^{i-1}I_{ji}I_{ik}\right]&\leq&\frac{n(n-1)(n-2)}{3}\\
&& (\mbox{here}\ \mbox{either}\ j<m=i<k,\ \mbox{or}\  m<k=j<i),
\end{eqnarray*}

\begin{eqnarray*}
E_3=2E\left[\sum_{i=3}^n\sum_{j=2}^{i-1}\sum_{m=1}^{j-1}I_{ji}I_{mi}\right]&\leq&\frac{n(n-1)(n-2)}{3}\ (\mbox{here}\ j,m<k=i,\ m\not= j),
\end{eqnarray*}

\begin{eqnarray*}
E_4=2E\left[\sum_{i=3}^n\sum_{k=2}^{i-1}\sum_{j=1}^{k-1}I_{ji}I_{jk}\right]&\leq&\frac{n(n-1)(n-2)}{3}\ (\mbox{here}\ k,i>m,\ m= j)
\end{eqnarray*}
and 
\begin{eqnarray*}
E_5&=&E\left[\sum_{i=2}^n\sum_{j=1}^{i-1}\ \sum_{k=2,\ k\not= j,i}^{n}\ \sum_{m=1,\ m\not=j,i}^{k-1}I_{ji}I_{mk}\right]\\
&=&\sum_{i=2}^n\sum_{j=1}^{i-1}\ \sum_{k=2,\ k\not= j,i}^{n}\ \sum_{m=1,\ m\not=j,i}^{k-1}P\left(C^{(n)}_{j,i;\ m,k}\right),
\end{eqnarray*}
where $1\leq j<i\leq n,\ 1\leq m<k\leq n,\ m,k\not=j,i$ and
$$
C^{(n)}_{j,i;\ m,k}=\{X_{j,n}\leq X_{i,n},\ Y_{j,n}\leq Y_{i,n},\ X_{m,n}\leq X_{k,n},\ Y_{m,n}\leq Y_{k,n}\}.
$$
Let
$$
D^{(n)}_{l,t;\ s,w}=\{X_t\leq X_l,\ Y_t\leq Y_l,\ X_w\leq X_s,\ Y_w\leq Y_s\}.
$$
We write as in the proof of Theorem\ref{theorem2.1} that
$$
\cup_{1\leq j<i\leq n,\ 1\leq m<k\leq n,\ m,k\not=j,i}\ C^{(n)}_{j,i;\ m,k}=\cup_{1\leq t,l,w,s\leq n,\ t\not= l,w,s,\ l\not= w,s,\ w\not= s}\ D^{(n)}_{l,t;\ s,w}.
$$
We have 
\begin{eqnarray*}
E_5&=&P\left(\cup_{1\leq t,l,w,s\leq n,\ t\not= l,w,s,\ l\not= w,s,\ w\not= s}\ D^{(n)}_{t,l;\ w,s}\right)\\
&=&\sum_{t=1}^n\ \sum_{l=1,\ l\not=t}^{n}\ \sum_{w=1,\ w\not= t,l}^{n}\ \sum_{s=1,\ s\not=t,l,w}^{n} p_{tl}p_{ws}.
\end{eqnarray*}
Let us study the last term in (\ref{2.1}). We have
\begin{eqnarray*}
\left(\sum_{k=1}^n\ \sum_{m=1,\ m\not= k}^n p_{mk}\right)^2&=&\sum_{k=1}^n\ \sum_{m=1,\ m\not= k}^n p_{mk}\sum_{t=1}^n\ \sum_{l=1,\ l\not= t}^n p_{tl}\\
&\geq& \sum_{k=1}^n\ \sum_{m=1,\ m\not= k}^n \sum_{t=1,\ t\not= m,k}^n\ \sum_{l=1,\ l\not= m,k,t}^n  p_{mk}p_{tl}=E_5.
\end{eqnarray*}
That way,
$$
Var(\tilde{\tau}_n)\leq \frac{16(E_1+E_2+E_3+E_4)}{(n-1)^4}\rightarrow 0\quad (n\rightarrow \infty).
$$
The result readily follows.
\end{gproof}

\section{Examples and Simulation Experiments}
{\bf 1.} Let 
$$
F_i(x,y)=\frac{1}{2\pi\sqrt{1-t_i^2}}\int_{-\infty}^x \int_{-\infty}^y e^{-\frac{u^2-2t_i uv+v^2}{2}}dudv\quad (x,y\in \mathbb{R},\ -1<t_i<1)
$$
be a bivariate normal distribution function with marginal distribution functions 
$$
H(x)=G(x)=\Phi(x)=\frac{1}{\sqrt{2\pi}}\int_{-\infty}^x e^{-\frac{u^2}{2}}du\quad (x\in \mathbb{R})
$$
and 
corresponding density functions $f_i(x,y)=F''_{i,xy}(x,y),\ \phi(x)=\Phi'(x)$. Let $F'_{i,x}(x,y)=\frac{\partial F_i(x,y)}{\partial x}$. Observe that
$$
F'_{i,x}(x,y)=\phi(x)\Phi\left(\frac{y-t_ix}{\sqrt{1-t^2_i}}\right),\quad F'_{i,y}(x,y)=\phi(y)\Phi\left(\frac{x-t_iy}{\sqrt{1-t^2_i}}\right).
$$
By integrating in parts by $y$, we get
\begin{eqnarray}\label{3.1}
\nonumber E(F_i(X_j,Y_j))&=& \int_{\mathbb{R}^2}F_i(x,y)f_j(x,y)dxdy\\
&=&1/2-\int_{\mathbb{R}^2} \phi(x) \phi(y) \Phi \left(\frac{x-t_iy}{\sqrt{1-t^2_i}}\right)\Phi \left(\frac{y-t_jx}{\sqrt{1-t^2_j}}\right)dxdy=Q(t_i,t_j).
\end{eqnarray}
Then
\begin{eqnarray*}
Q'_{t_i}(t_i,t_j)&=&\int_{\mathbb{R}^2} \phi(x) \phi(y) \phi \left(\frac{x-t_iy}{\sqrt{1-t^2_i}}\right)\left(\frac{y-t_ix}{(1-t^2_i)^{3/2}}\right)\Phi \left(\frac{y-t_jx}{\sqrt{1-t^2_j}}\right)dxdy\\
&=&\int_{\mathbb{R}^2} (\phi(x) )^2 \phi \left(\frac{y-t_ix}{\sqrt{1-t^2_i}}\right)\left(\frac{y-t_ix}{(1-t^2_i)^{3/2}}\right)\Phi \left(\frac{y-t_jx}{\sqrt{1-t^2_j}}\right)dxdy\\
&=&-\frac{1}{\sqrt{1-t_i^2}}\int_{\mathbb{R}} (\phi(x) )^2 \left(\int_{\mathbb{R}} \Phi \left(\frac{y-t_jx}{\sqrt{1-t^2_j}}\right) d_y \phi \left(\frac{y-t_ix}{\sqrt{1-t^2_i}}\right)\right) dx.\\
\end{eqnarray*}
By integrating in parts again, we obtain
\begin{eqnarray*}
Q'_{t_i}(t_i,t_j)&=&\frac{1}{\sqrt{1-t_i^2}\sqrt{1-t_j^2}}\int_{\mathbb{R}} (\phi(x) )^2 \left(\int_{\mathbb{R}} \phi \left(\frac{y-t_jx}{\sqrt{1-t^2_j}}\right) \phi \left(\frac{y-t_ix}{\sqrt{1-t^2_i}}\right)dy \right) dx.\\
\end{eqnarray*}
Observe that
\begin{eqnarray*}
\int_{\mathbb{R}} \phi \left(\frac{y-t_jx}{\sqrt{1-t^2_j}}\right)&\cdot& \phi \left(\frac{y-t_ix}{\sqrt{1-t^2_i}}\right)dy\\
&=&\frac{e^{x^2}\sqrt{1-t_i^2}\sqrt{1-t_j^2}}{\sqrt{2\pi}\sqrt{2-t_i^2-t_j^2}}\cdot e^{-\frac {x^2[(2-t_i^2-t_j^2)^2-(t_i(1-t_j^2)+t_j(1-t_i)^2)^2]}  {2(1-t_i^2)(1-t_j^2)(2-t_i^2-t_j^2)}}.
\end{eqnarray*}
Then
$$
Q'_{t_i}(t_i,t_j)=\frac{\sqrt{1-t_i^2}\sqrt{1-t_j^2}}{2\pi\sqrt{(2-t_i^2-t_j^2)^2-(t_i(1-t_j^2)+t_j(1-t_i^2))^2}}.
$$
After some algebra, we obtain 
\begin{equation}\label{3.2}
Q'_{t_i}(t_i,t_j)=\frac{1}{2\pi\sqrt{1-(t_i+t_j)^2}}.
\end{equation}
Putting $t_i=t_j$ in (\ref{3.1}), we get
\begin{eqnarray*}
Q(t_j,t_j)=1/2-\int_{\mathbb{R}^2} \phi(x) \phi(y) \Phi \left(\frac{x-t_jy}{\sqrt{1-t^2_j}}\right)\Phi \left(\frac{y-t_jx}{\sqrt{1-t^2_j}}\right)dxdy.
\end{eqnarray*}
It was shown, for example, in Stepanov (2025a) that
$$
Q(t_j,t_j)=\frac{\arcsin (t_j)}{2\pi}+1/4.
$$
By integrating (\ref{3.2}), we obtain that
$$
Q(t_i,t_j)=\frac{\arcsin \left(\frac{t_i+t_j}{2}\right)}{2\pi}+C.
$$
It follows from the last two equalities that  $C=1/4$. We finally get
$$
EF_i(X_j,Y_j)=\frac{\arcsin \left(\frac{t_i+t_j}{2}\right)}{2\pi}+1/4.
$$
Then
\begin{eqnarray*}
\tau_n&=&\frac{2}{\pi n(n-1)}\sum_{i=1}^n\ \sum_{j=1,\ j\not= i}^n \arcsin \left(\frac{t_i+t_j}{2}\right)\\
&=& \frac{4}{\pi n(n-1)}\sum_{i=2}^n\ \sum_{j=1}^{i-1} \arcsin \left(\frac{t_i+t_j}{2}\right).
\end{eqnarray*}
Observe that when $t_i=t_j=t$, we get $\tau_n=\tau= \frac{2}{\pi}\arcsin (t)$, see, for example Stepanov(2025a).

\noindent {\bf 1.1}\  Let $t_i=\sin(i)\ (i\geq 1)$. For $n=100000$ we have found that $\tau_n=2.4438*10^{-7}$. Obviously $\tau=\lim_{n\rightarrow \infty}\tau_n=0$. For the same $n$ we have conducted a simulation experiment. We have obtained that $\tilde{\tau}_n=-0.0013$.

\noindent {\bf 1.2}\  Let now $t_i=e^{-\mid \sin(i)\mid}\ (i\geq 1)$. For $n=100000$ we have found that $\tau_n=0.3826$ and $\tilde{\tau}_n=0.3882$.

\noindent {\bf 2.}  Let us consider a Farlie-Gumbel-Morgenstern copula with distribution and density functions
$$
F_i(x,y)=xy+t_i(x-x^2)(y-y^2),
$$$$
f_i(x,y)=1+t_i(1-2x)(1-2y),
$$
where $0<x,y<1,\ -1\leq t_i\leq 1$. We have
$$
E(F_i(X_j,Y_j))=\frac{n(n-1)}{4}+\frac{(n-1)\sum_{i=1}^n t_i}{18}.
$$
Then
$$
\tau_n=\frac{2\sum_{i=1}^n t_i}{9n}.
$$
\noindent {\bf 2.1}\ Let $t_i=1/i$. Then $\tau=0$. We have conducted a simulation experiment for $n=100000$ and obtained $\tilde{\tau}_n=-0.0022$.

\noindent {\bf 2.2} Let $t_i=3/5-1/i$. Then $\tau=2/15\approx 0.1333$. We have conducted  a simulation experiment for $n=100000$ and obtained $\tilde{\tau}_n=0.1317$.

\noindent {\bf 3.} Let now $F_i$ be a bivariate Pareto distribution function with
$$
F_i(x,y)=1-\frac{1}{(1+x)^{t_i}}-\frac{1}{(1+y)^{t_i}}+\frac{1}{(1+x+y)^{t_i}}
$$
and
$$ 
f_i(x,y)=\frac{t_i(t_i+1)}{(1+x+y)^{t_i+2}},
$$
where $x,y,t_i>0$. We have
$$
E(F_i(X_j,Y_j))=\frac{t_j^2+t_j}{(t_i+t_j)(t_i+t_j+1)}
$$
and
$$
\tau_n=\frac{4\sum_{j=1}^n(t_j^2+t_j)\sum_{i=1,\ i\not= j}^n \frac{1}{(t_i+t_j)(t_i+t_j+1)}}{n(n-1)}-1.
$$
Let $t_i=i$. Then 
$$
\tau_n=\frac{4}{n(n-1)}\sum_{j=1}^n\frac{j(j+1)(n-j)}{(2j+1)(n+j+1)}\rightarrow 0.2275\quad (n\rightarrow \infty).
$$
We have conducted a simulation experiment with $n=100000$. We obtained $\tilde{\tau}_n=0.2297$.

It should be noted that all simulation experiments support our theoretical results.

%Disclosure statement

%No potential conflict of interest was reported by the authors.

\section*{References}
\begin{description} 

\item Bairamov, I., Stepanov, A. (2010).\ Numbers of near-maxima for the bivariate case,  {\it Statistics  $\&$ Probability Letters}, {\bf 80}, 196--205.

\item Balakrishnan, N., Lai, C.D. (2009).\ {\it Continuous Bivariate Distributions}, Second edition, Springer.
 
\item Balakrishnan, N., Stepanov, A. (2015).\ Limit results for concomitants of order statistics, {\it Metrika}, {\bf 78}, 385--397.

\item  Bhattacharya, B.B. (1974).\ Convergence of sample paths of normalized sums of induced order statistics, {\it Ann. Statist. }, {\bf 2}, 1034--1039.

\item David, H.A., Galambos, J. (1974).\ The asymptotic theory of concomitants of order statistics, {\it J. Appl. Probab.}, {\bf 11}, 762--770.

\item David, H.A., Nagaraja, H.N. (2003).\  {\it Order Statistics}, Third  edition, John Wiley \& Sons, NY.

\item Egorov, V. A., Nevzorov, V. B. (1984).\ Rate of convergence to the Normal law of sums of induced order statistics,
{\it Journal of Soviet Mathematics} (New York), {\bf 25}, 1139--1146.

\item Fisher, R.A. (1921).\ On the ''probable error" of a coefficient of correlation deduced from a small sample, Metron, {\bf 1},  3--32.

\item Hauke, J.,  Kossowski, T. (2011).\ Comparison of values of Pearson's and Spearman's correlation coefficients on the same sets of data, {\it Quaestiones Geographicae}, {\bf 30} (2), 87--93.

\item Hung, Y. C., Chen, R., Balakrishnan, N. (2016).\  On the correlation structure of exponential order statistics and some extensions, {\it Mathematical Methods of Statistics}, {\bf 25},  196--206.

\item Kendall, M. G. (1970).\ {\it Rank Correlation Methods}, London, Griffin.

\item Nevzorov, V. B. (2001). {\it Records: Mathematical Theory}, Translation of Mathematical Monographs, Vol. {\bf 194}, American Mathematical Society, Providence, Rhode Island.

\item Shevlyakov, G.L., Vilchevski, N.O. (2002).\ {\it Robustness in Data Analysis: Criteria and Methods, Modern Probability and Statistics}, VSP, Utrecht.

\item Stepanov, A.V. (2025a).\ On Kendall's correlation coefficient, {\it Vestnik St. Petersburg University, Mathematics}, {\bf 58} (1), 71--78.

\item Stepanov, A. (2025b).\ Comparison of Correlation Coefficients,  {\it Sankhya A}, {\bf 87} (1),  191--218.

\item  Stepanov, A. (2025c).\ On Rank Correlation Coefficients, {\it arXiv:2506.06056v2}.

\item  Xu, W., Hou,  Y.,   Hung, Y. S., Zou, Y. (2013).\ A comparative analysis of Spearman's rho and Kendall's tau
in normal and contaminated normal models, {\it Signal Processing}, {\bf 93}, 261–-276.

\end{description}

\end{document}